\documentclass[conference]{IEEEtran} 

\usepackage{amsmath}
\usepackage{amssymb}
\usepackage{amsthm}
\usepackage{latexsym}
\usepackage{MnSymbol}
\usepackage{verbatim}
\usepackage[final]{graphicx}
\usepackage{psfrag}
\usepackage{algorithm}
\usepackage{algorithmic}

\newtheorem{proposition}{Proposition}

\newtheorem{remark}{Remark}

{\begin{list}               
    {$\bullet$ \hfill}{
        \setlength{\leftmargin}{\parindent}
        \setlength{\parsep}{0.04\baselineskip}
        \setlength{\itemsep}{0.5\parsep}
        \setlength{\labelwidth}{\leftmargin}
        \setlength{\labelsep}{0em}}
    }
{\end{list}}

\providecommand{\eref}[1]{\eqref{eq:#1}}  
\providecommand{\cref}[1]{Chapter~\ref{chap:#1}}
\providecommand{\sref}[1]{Section~\ref{sec:#1}}

\providecommand{\R}{\ensuremath{\mathbb{R}}}

\providecommand{\norm}[1]{\lVert#1\rVert}

\providecommand{\bydef}{\overset{\text{def}}{=}}
\providecommand{\tr}{\operatorname{tr}}
\providecommand{\opvec}{\operatorname{vec}}
\providecommand{\opmat}{\operatorname{mat}}

\renewcommand{\vec}[1]{\ensuremath{\boldsymbol{#1}}}
\providecommand{\mat}[1]{\ensuremath{\boldsymbol{#1}}}

\providecommand{\mA}{\mat{A}} \providecommand{\mB}{\mat{B}}
\providecommand{\mC}{\mat{C}} \providecommand{\mH}{\mat{H}}
\providecommand{\mD}{\mat{D}} 
\providecommand{\mI}{\mat{I}}  
  
\providecommand{\mM}{\mat{M}} \providecommand{\mP}{\mat{P}} 
\providecommand{\mQ}{\mat{Q}}

\providecommand{\va}{\vec{a}} 
 \providecommand{\ve}{\vec{e}}

\providecommand{\vx}{\vec{x}} \providecommand{\vy}{\vec{y}}
\providecommand{\vz}{\vec{z}} 
 \providecommand{\vzero}{\vec{0}}
\providecommand{\veta}{\vec{\eta}}
 \providecommand{\vv}{\vec{v}}

\providecommand{\ve}{\vec{e}}
\providecommand{\vf}{\vec{f}}

\providecommand{\projp}[1]{\ensuremath{\mP^{\perp}_{#1}}}

\providecommand{\lmax}{\lambda_{\max}}

\usepackage[nosort]{cite}
\usepackage{algorithmic}
\usepackage{algorithm}

\usepackage{pgfplots}
\usepackage{tikz}

\usepackage{subfigure}
\usepackage{caption}
\newcommand\hairspace{\kern .08333em }
 
\providecommand{\QQ}{\mQ}
\providecommand{\PP}{\mP}
\providecommand{\DD}{\mD}
\providecommand{\ee}{\vec{e}}
\providecommand{\ff}{\vec{f}}
\providecommand{\MSE}[1]{\ensuremath{\E \left|\left|\vz^{(#1)}\right|\right|^2}}

\providecommand{\lmax}{\lambda_{\max}}
\providecommand{\E}{\mathbb{E} \hairspace}
\providecommand{\En}{\mathbb{E}_{\veta} \hairspace}
\providecommand{\projp}[1]{\ensuremath{\mP^{\perp}_{#1}}}
\providecommand{\tva}{\widetilde{\va}}
\providecommand{\teta}{\widetilde{\eta}}

\providecommand{\vzs}[1]{[\vz^{(#1)}]^{\otimes 2}}
\providecommand{\HH}{\begin{pmatrix} \QQ & \DD \\ \mat{0} & \PP \end{pmatrix}}

\title{Randomized Kaczmarz Algorithm for Inconsistent Linear Systems: An Exact MSE Analysis}
\IEEEoverridecommandlockouts

\author{
\IEEEauthorblockN{Chuang Wang$^1$, Ameya Agaskar$^{1,2}$~and~Yue M.~Lu$^1$\thanks{The Lincoln Laboratory portion of this work was sponsored by the Department 
                         of the Air Force under Air Force Contract \#FA8721-05-C-0002. 
                         Opinions, interpretations, conclusions and recommendations are those of the 
                         authors and are not necessarily endorsed by the United States Government.
                          }%
                          \thanks{C. Wang and Y.~M.~Lu were supported in part by the U.S. National Science Foundation under Grant CCF-1319140.}
}
\vspace{1ex}
\IEEEauthorblockA{$^1$Harvard University, Cambridge, MA 02138, USA \\
                  $^2$MIT Lincoln Laboratory, Lexington, MA 02420, USA \\
                  E-mail:\{chuangwang,aagaskar,yuelu\}@seas.harvard.edu 
                  }
}

\begin{document}
\bstctlcite{NOURL}

%
\maketitle
\begin{abstract}
We provide a complete characterization of the randomized Kaczmarz algorithm (RKA) for inconsistent linear systems.
The Kaczmarz algorithm, known in some fields as the algebraic reconstruction technique,
is a classical method for solving large-scale overdetermined linear systems through a sequence of projection
operators; the randomized Kaczmarz algorithm is a recent proposal by Strohmer and Vershynin to randomize
the sequence of projections in order to guarantee exponential convergence (in mean square) to the solutions. A flurry of work
followed this development, with renewed interest in the algorithm, its extensions, and various bounds
on their performance. Earlier, we studied the special case of consistent linear systems and 
provided an exact formula for the mean squared error (MSE) in the value reconstructed
by RKA, as well as a simple way to compute the exact decay rate of the error. In this work, we consider
the case of inconsistent linear systems, which is a more relevant scenario for most applications. First, by using a ``lifting trick'', we derive an exact formula for the MSE given a fixed noise vector added to the measurements. Then we show how
to average over the noise when it is drawn from a distribution with known first and second-order statistics. Finally, we demonstrate the accuracy of our exact MSE formulas through numerical simulations, which also illustrate that previous upper bounds in the literature may be several orders of magnitude too high.
\end{abstract}

\section{Introduction}
The Kaczmarz algorithm \cite{Kaczmarz:37} is a simple and popular iterative method for solving large-scale overdetermined linear systems.
Given a full-rank measurement matrix $\mA \in \R^{m \times n}$, with $m \geq n$, we wish to recover a signal vector $\vx \in \R^n$ from its measurements $\vy \in \R^m$, given by
\begin{equation}
 \vy = \mA \vx.
\label{eq:yeqAx}
\end{equation}
Each row of $\mA$ describes a single linear measurement $y_i = \va_i^T \vx$, and the set of all signals $\vx$ that satisfy that equation is an $(n-1)$-dimensional (affine) subspace in $\R^n$.
The Kaczmarz algorithm begins with an arbitrary initial guess $\vx^{(0)}$, then cycles through all the rows, projecting the iterand $\vx^{(k-1)}$ onto the subspace
$\{\vx \in \R^n : \va_r^T \vx = y_r\}$ to obtain $\vx^{(k)}$, where $r = k \operatorname{mod} m$. Since affine subspaces are convex, this algorithm is a special case of the projection onto convex sets (POCS) algorithm \cite{trussell_signal_1984}.

Due to its simplicity, the Kaczmarz algorithm has been widely used in signal and image processing
. It has long been observed by practitioners that its convergence rate depends on the ordering of the rows
of \mA, and that choosing the row order at random can often lead to faster convergence \cite{herman_algebraic_1993}.
Yet it was only recently that Strohmer and Vershynin first rigorously analyzed the \textit{randomized} Kaczmarz algorithm (RKA) \cite{strohmer_randomized_2009}. They considered a convenient row-selection probability distribution: choosing row $i$ with probability proportional to its squared norm $||\va_i||^2$, and proved the following upper bound on the mean squared error (MSE) of the RKA at the $k$th iteration:
\begin{equation}\label{eq:strohmer_bounds}
\mathbb{E} \, \norm{\vx^{(k)} - \vx}^2 \leq \left(1 - \kappa_{\mA}^{-2}\right)^k \norm{\vx^{(0)} - \vx}^2,
\end{equation}
where $\kappa_{\mA} \bydef \norm{\mA}_F \norm{\mA^{-1}}_2$ is related to the condition number of $\mA$, and $\mA^{-1}$ is its left-inverse. 
This bound guarantees that the MSE decays \emph{exponentially} as the RKA iterations proceed.

The work of Strohmer and Vershynin spurred a great deal of interest in RKA and its various extensions (see,~\emph{e.g.}, 
\cite{censor_note_2009, needell_randomized_2010,Chen:2012fe, zouzias_randomized_2013,needell_paved_2014, Dai:2014, recht_toward_2012}). 
In particular, Needell \cite{needell_randomized_2010} derived a bound for the performance of the algorithm when the underlying linear system is inconsistent. In this case, the measurements are 
\begin{equation}
\vy = \mA \vx + \veta,
\label{eq:yeqAxpluseta}
\end{equation}
where $\veta$ is an additive noise vector. Needell's bound was later improved by Zouzias and Freris \cite{zouzias_randomized_2013}, who proved the following upper bound on the MSE:
\begin{equation}
    \mathbb{E} \, \norm{\vx^{(k)} - \vx}^2 \leq \left(1 - \kappa_{\mA}^{-2}\right)^k \norm{\vx^{(0)} - \vx}^2 + \frac{||\veta||^2}{\sigma^2_{\min}(\mA)},
\label{eq:zf_bound}
\end{equation}
where $\sigma_{\min}(\mA)$ is the smallest singular value of $\mA$. Note that this bound is equal to the original noiseless bound in (\ref{eq:strohmer_bounds}) plus an extra term proportional to the total squared error in the measurements.

In this paper, we provide a complete characterization of the randomized Kaczmarz algorithm for inconsistent linear systems given in \eref{yeqAxpluseta}.  We show in \sref{fixed} how to compute the \emph{exact} MSE of the algorithm (averaging over the random choices of the rows) at each iteration. This extends our earlier results derived for the special case when the measurements are noiseless \cite{agaskar_randomized_2014}. A key ingredient of our derivation is a ``lifting trick'', which allows us to analyze the evolution of the MSE (a quadratic quantity) through a much simpler linear recursion embedded in a higher-dimensional ``lifted'' space. We show that existing upper bounds in the literature \cite{needell_randomized_2010, zouzias_randomized_2013} can be easily derived from our exact MSE formula. By setting the number of iterations to infinity, we also provide a closed-form expression for the limiting MSE (\emph{i.e.}, error floor) of the algorithm.

Our MSE analysis in \sref{fixed} is conditioned on the noise, \emph{i.e.}, we assume that $\veta$ in \eref{yeqAxpluseta} is a \emph{fixed} and \emph{deterministic} vector. Thus, the resulting expressions for the MSE and the error floor depend on $\veta$. In practice, the measurement noise $\veta$ is unknown but its elements can often be modeled as zero-mean i.i.d. random variables drawn from some probability distributions (\emph{e.g.}, Gaussian random variables.) We consider this setting in Section \ref{sec:gauss}, where we compute the MSE with the expectations taken over two sources of randomness: the random row-selections made by the algorithm and the noise vector. In this case the final MSE has two terms: one equals to that of the noiseless case we analyzed in\cite{agaskar_randomized_2014} and an extra term proportional to the noise variance. 

We demonstrate the accuracy of our exact MSE formulas through numerical simulations reported in \sref{results}. These empirical results also illustrate that previous upper bounds in the literature may be several orders of magnitude too high over the true performance of the algorithm.


\section{Exact Performance Analysis of RKA}
\label{sec:fixed}

\subsection{Overview and Notation}

Consider an inconsistent linear system as in \eref{yeqAxpluseta}. Given the measurement $\vy$ and the matrix $\mA$, the randomized Kaczmarz algorithm seeks to (approximately) reconstruct the unknown signal $\vx$ through iterative projections. The iterand $\vx^{(0)} \in \R^n$ is initialized arbitrarily. At the $k$th step, a row $i_k$ is chosen at random; row $i$ is chosen with probability $p_i$. The update equation is
\begin{equation}
\vx^{(k)} = \vx^{(k-1)} + \frac{y_{i_k} - \va_{i_k}^T\vx^{(k-1)}}{||\va_{i_k}||^2}\va_{i_k},
\label{eq:kacz-itr}
\end{equation}
where $\va_i$ is the $i$th row of $\mA$. The intuition behind the algorithm is simple. In the noiseless case (\emph{i.e.}, when $\veta = \vzero$), each row of $\mA$ and its corresponding entry in $\vy$ defines an affine subspace on which the solution $\vx$ must lie; at each iteration, the RKA algorithm randomly selects one of these subspaces and projects the iterand onto it, getting closer to the true solution with each step. 

The row-selection probabilities $p_i$ are tunable parameters of the algorithm. Other authors \cite{strohmer_randomized_2009, needell_randomized_2010, zouzias_randomized_2013} have fixed the probabilities to be $p_i = ||\va_i||^2/||\mA||_F^2$. This is not really a restriction, since the rows of $\mA$ can be scaled arbitrarily and, as long as the measurements are scaled appropriately, the solution and the algorithm's iterations do not change. This particular choice, though, is convenient because it leads to simplified bounds, allowing them to be written in terms of a (modified) condition number of $\mA$. Since most of our expressions are not simplified through this choice, we fix $\mA$ and allow the $p_i$ to be chosen arbitrarily. Our only restriction is that $p_i > 0$ for all $i$ so that every measurement is used.


In \cite{agaskar_randomized_2014}, we computed the exact MSE of the RKA for the special case of consistent linear systems (\emph{i.e.}, the noiseless case.) In what follows, we extend our earlier result and analyze the more general inconsistent case.

%
%
%


\begin{table}
    \renewcommand{\arraystretch}{1.8}
    \centering
    \begin{tabular}{l|l||l|l}
	Variable & Definition & Variable & Definition\\
        \hline
        $\tva$          & $\frac{\va}{ \|\va \|}$ &
        $\projp{i}$ & $\mI - \tva_i \tva_i^T$ 
        \\
        $\teta_i$       & $\frac{\eta_i}{\| \va_i \|}$ &
        $\vv^{\otimes 2}$ & $\vv \otimes \vv$
        \\
        $\PP$           & $\E  \projp{i}$ &
        $\QQ$           & $\E  \projp{i} \otimes \projp{i}$ 
        \\
        $\vf$           & $\E  \tva_{i} \teta_{i}  $ &
        $\ve$           & $\E   \teta_i^{\hairspace 2} \tva_i \otimes \tva_i$ 
        \\
        $\DD$           & $\E  \mD_i \teta_i $&
        $\mD_i$         & $\tva_i \otimes \projp{i}
                           + \projp{i} \otimes \tva_i $
        \\
        $\vv_2$         & $(\mI - \PP)^{-1} \vf$ &
        $\vv_1$         & $(\mI - \QQ)^{-1} \left[ \ve + \DD \vv_2 \right]$
    \end{tabular}
    \caption{List of important notation and variables.}
    \label{tb:var}
\end{table}

\subsection{Exact MSE Analysis Using Lifting}

We consider a given measurement matrix $\mA$ and a fixed noise vector $\veta$.
The exact MSE at iteration step $k$ over the randomness of the algorithm are formulated in this section. 

To lighten the notation, we define the normalized $i$th row vector of $\mA$ 
as $\tva_i \bydef \frac{\va_i}{\| \va_i \|} \in \R^n$ and let $\teta_i \bydef \frac{\eta_i}{\| \va_i \|}$.
These and other important definitions are summarized in Table \ref{tb:var} for the the reader's convenience.
By combining \eqref{eq:kacz-itr} and \eqref{eq:yeqAxpluseta}, 
the error vector $\vz^{(k)} = \vx^{(k)} - \vx$ can be expressed as
\begin{equation}
  \vz^{(k)} = \projp{i_k} \vz^{(k-1)} + \tva_{i_k} \teta_{i_k},
  \label{eq:z}
\end{equation}
where $\projp{i} \bydef \mI - \tva_i \tva_i^T$ is the projection onto the $(n-1)$-dimensional subspace orthogonal to the $i$th row $\tva_i$.
Averaging \eqref{eq:z} over the randomness of the algorithm, we get an iterative equation of the mean error vector
\begin{equation}
  \E  \vz^{(k)} = \PP \, \E  \vz^{(k-1)} + \vf,
  \label{eq:ez}
\end{equation}
where $\PP \bydef \E  \projp{\va_{i}}$ and  $\vf \bydef \E\hairspace  \tva_{i} \teta_{i}$.

Note that the ease with which we can obtain \eref{ez} from \eref{z} is mainly due to the \emph{linearity} of the original random recursion in \eref{z}. However, the quantity we are interested in, the mean-squared error $\MSE{k}$, is a non-linear (quadratic) term. To compute it, we ``lift'' the problem by treating the covariance matrix $\E \vz^{(k)} (\vz^{(k)})^T$ as an $n^2$-dimensional vector whose dynamics are determined by the algorithm. In the lifted space, the dynamics
are still linear (as we shall soon see), thus allowing for a relatively simple analysis.  The MSE can be easily obtained as the trace of $\E \vz^{(k)} (\vz^{(k)})^T$.

Consider the $k$th iteration:
\begin{equation}
\begin{aligned}\label{eq:zzT}
    &\vz^{(k)} (\vz^{(k)})^T  = \projp{i_k} \vz^{(k-1)} (\vz^{(k-1)})^T \projp{i_k} + \teta_{i_k}^2 \tva_{i_k} \tva_{i_k}^T\\
       & \qquad\qquad + \teta_{i_k} \left(\tva_{i_k} (\vz^{(k-1)})^T \projp{i_k} + \projp{i_k} \vz^{(k-1)} \tva_{i_k}^T\right).
\end{aligned}
\end{equation}
The linearity of this expression will be clearer if we ``vectorize'' $\vz^{(k)} (\vz^{(k)})^T$ by vertically
concatenating its columns to form a vector $ \opvec\left(\vz^{(k)} (\vz^{(k)})^T\right) \in \R^{n^2}$. In what follows, we will make use of the following matrix identity which holds for any matrices dimensioned so that $\mA \mB \mC$ is well-defined:
\begin{equation}
    \opvec(\mA \mB \mC) = (\mC^T \otimes \mA) \opvec(\mB),
    \label{eq:veckronidentity}
\end{equation}
where $\otimes$ represents the Kronecker matrix product \cite{horn_topics_1994}. First, we note that
\[ 
\opvec\left( \vz^{(k)} (\vz^{(k)})^T\right) = \vz^{(k)} \otimes \vz^{(k)} \bydef \vzs{k}, 
\]
where $\vv^{\otimes 2}$ is introduced as a shorthand notation for the Kronecker product of a vector $\vv$ and itself. Then, we can apply the identity \eqref{eq:veckronidentity} to the
right hand side of \eqref{eq:zzT} to obtain
\[
\begin{aligned}
    \vzs{k} & = \left(\projp{i_k} \otimes \projp{i_k}\right) \vzs{k-1} +  \teta^2_{i_k} \tva_{i_k}^{\otimes 2}\\
    & \qquad +  \teta_{i_k} \left( \projp{i_k} \otimes \tva_{i_k} + \tva_{i_k} \otimes \projp{i_k}  \right) \vz^{(k-1)}.
\end{aligned}
\]
Taking expectation on both sides of the equation over the randomness of the algorithm, we obtain a simple iterative
formula for the second-moment matrix:
\begin{equation}
    \E \vzs{k} = \QQ \, \E \vzs{k-1} + \DD \, \E \vz^{(k-1)} + \ve, 
    \label{eq:ezz}
\end{equation}
where $\QQ \bydef \E  \left(\projp{i} \otimes \projp{i}\right)$,
$\DD \bydef \E \teta_{i} \left( \projp{i} \otimes \tva_{i} + \tva_{i} \otimes \projp{i}  \right)$
and $\ve \bydef \E \teta^{\hairspace 2}_{i} \, \tva_{i}^{\otimes 2}$.

We can combine \eqref{eq:ez} and \eqref{eq:ezz} into a single linear recursion
\begin{equation}\label{eq:HH}
    \begin{pmatrix}
      \E \vzs{k}\\
      \E \vz^{(k)}
    \end{pmatrix}
     =
    \mH
    \begin{pmatrix}
      \E \vzs{k-1} \\
      \E \vz^{(k-1)}
    \end{pmatrix}
    +
    \begin{pmatrix}
      \ee\\
      \ff
    \end{pmatrix},
\end{equation}
where 
\begin{equation}\label{eq:H_def}
\mH \bydef \HH.
\end{equation}

We thus have the following proposition:
\begin{proposition}\label{prop:MSE}
  For a fixed noise vector $\veta$, and an initial error vector $\vz^{(0)}$, the MSE of RKA at the $k$th iteration is given by
\begin{equation}\label{eq:fixed_noise_iteration_MSE}
\MSE{k}= \begin{pmatrix}
    \opvec(\mI_{n})\\
    \vec{0}_{n}
  \end{pmatrix}^{T}
  \left[
\mH^{k}
    \begin{pmatrix}
      \vzs{0} - \vv_{1}\\
      \vz(0) - \vv_{2}
    \end{pmatrix}
    +
    \begin{pmatrix}
      \vv_1\\
      \vv_2
    \end{pmatrix}
  \right],
\end{equation}
where 
$
  \vv_2 = (\mI - \PP)^{-1} \vf
$ and
$
  \vv_1 = (\mI - \QQ)^{-1} \left[ \ee + \DD \vv_2  \right]
$.
\end{proposition}
\begin{IEEEproof}
We first solve the linear recursion \eref{HH} to get a closed-form expression
\begin{align}
     \begin{pmatrix}
      \E \vzs{k}\\
      \E \vz^{(k)}
    \end{pmatrix}
     & =
    \mH^k
    \begin{pmatrix}
      \E \vzs{0} \\
      \E \vz^{(0)}
    \end{pmatrix}
    +
    \sum_{\ell=0}^{k-1}
    \mH^\ell
    \begin{pmatrix}
      \ee \\
      \ff
    \end{pmatrix}
    \label{eq:solved_H_iteration}
\end{align}
that depends on the initial error $\vz^{(0)}$. Using the identity
$
    \sum_{\ell=0}^{k-1}
    \mH^\ell
    =(\mI-\mH^k)(\mI-\mH)^{-1}
$ and noting that 
\[
    (\mI - \mH)^{-1} = \begin{pmatrix}
      (\mI - \QQ )^{-1} & (\mI - \QQ)^{-1} \DD  (\mI - \PP)^{-1} \\
      \mat{0} & (\mI - \PP )^{-1}
    \end{pmatrix},
\]
we can simplify \eref{solved_H_iteration} and get
\begin{equation}\label{eq:solved_H2}
    \begin{pmatrix}
      \E \vzs{k}\\
      \E \vz^{(k)}
    \end{pmatrix}
     = 
\mH^{k}
    \begin{pmatrix}
      \vzs{0} - \vv_{1}\\
      \vz^{(0)} - \vv_{2}
    \end{pmatrix}
    +
    \begin{pmatrix}
      \vv_1\\
      \vv_2
    \end{pmatrix}.
\end{equation}

Meanwhile, using \eqref{eq:veckronidentity} and the fact that $(A \otimes B)^T = (A^T \otimes B^T)$, the MSE can be expressed 
in terms of the vectorized second-moment matrix as 
\begin{equation}
  \MSE{k} = \opvec({\mI _n})^T \left( \E \vzs{k} \right).
  \label{eq:mse-ezz}
\end{equation}
Combining this with \eref{solved_H2} yields the desired result.
\end{IEEEproof}

\subsection{The Limiting MSE}

With the iteration number $k$ going to infinity, the MSE in \eref{fixed_noise_iteration_MSE} will converge to a limiting value (\emph{i.e.}, an error floor) that only depends on the error vector $\veta$. To see this, we first note that both $\mP$ and $\mQ$ are positive semidefinite matrices by their constructions. In fact, one can show that
\[
0 < \lambda_{\max}(\QQ) \le \lambda_{\max}(\mP) < 1,
\]
where $\lambda_{\max}(\cdot)$ is the largest eigenvalue of matrix. Furthermore, one can show that the set of eigenvalues of $\mH $ defined in \eref{H_def} is the union of those of $\mP$ and $\mQ$, and that $\mH$ is a contraction mapping, with $\lim_{k \rightarrow \infty} \mH^k = \vzero$.  

This contraction property of $\mH$ implies that the first term in the right-hand side of \eref{fixed_noise_iteration_MSE} vanishes as $k$ goes to infinity. It follows that the limiting MSE can be characterized through $\vv_1$ as follows.

\begin{proposition}
The limiting MSE is given by 
  \begin{align}
    \MSE{\infty} = \tr \opmat(\vv_1),
    \label{eq:mse-infty}
  \end{align}
  where the $\opmat$ operator undoes the $\opvec$ operator to produce an $n \times n$ matrix from an $n^2$-dimensional vector.
\end{proposition}

Due to the space limit, we omit the proof of these assertions, which involve elementary matrix analysis and will be presented in a follow-up paper.


\subsection{Previous Upper Bounds on the MSE}
\label{subsec:previousbounds}
The existing bounds on noisy Kaczmarz performance \cite{needell_randomized_2010,zouzias_randomized_2013} can be recovered via our formulation. 
From \eqref{eq:HH}, we have
\begin{equation}\label{eq:recursion_MSE}
\begin{aligned}
    \MSE{k} & = \opvec(\mI)^T \QQ \opvec\left(\E \vz^{(k-1)} (\vz^{(k-1)})^T\right) \\
            & \quad\qquad + \opvec(\mI)^T \DD \, \E \vz^{(k-1)} + \opvec(\mI)^T \ee.
\end{aligned}
\end{equation}
Using the definition of $\mQ$, we have
\begin{align}
    \opvec(\mI)^T \QQ &=\left(\sum_i p_i (\projp{i} \otimes \projp{i}) \opvec(\mI)\right)^T\nonumber\\
&=\opvec(\sum_i p_i \projp{i} \projp{i})^T\nonumber\\
&=\opvec(\mP)^T,\nonumber
\end{align}
where the second equality can be obtained from the identity \eref{veckronidentity} and the last equality follows from $\projp{i}$ being an idempotent matrix and from the definition of $\mP$. It follows that the first term on the right-hand side of \eref{recursion_MSE} can be bounded as follows:
\[
\begin{aligned}
&\opvec(\mI)^T \QQ \opvec\left(\E \vz^{(k-1)} \vz^{(k-1)T}\right)\\
&\quad\qquad  =\opvec(\mP)^T \opvec\left(\E \vz^{(k-1)} \vz^{(k-1)T}\right)\nonumber \\
& \quad\qquad = \E \vz^{(k-1)T} \PP \vz^{(k-1)} \leq \lambda_{\max}(\PP) \E ||\vz^{(k-1)}||^2. \nonumber
\end{aligned}
\]

The second term on the right-hand side of \eref{recursion_MSE} is $0$, since
\[
    \DD^T \opvec(\mI) = \sum_i p_i \teta_i \left[\opvec(\projp{i} \tva_i) + \opvec(\tva_i^T \projp{i}) \right] = \vec{0}.
\]

The third term is given by
\[
    \opvec(\mI)^T \ve = \opvec(\mI)^T \sum_i p_i \teta_i^{\hairspace 2} \hairspace \tva_i^{\otimes 2} = \sum_i p_i \teta_i^{\hairspace 2}.
\]

So, all together, we have
\[
    \MSE{k} \leq \lambda_{\max}(\mP) \MSE{k-1} + \sum_i p_i \teta_i^2.
\]
Applying this inequality recursively gives us a bound equivalent to that in Zouzias and Freris \cite{zouzias_randomized_2013}:
\begin{align}
  \MSE{k}
   \leq \lmax^k(\PP) \MSE{0} + \frac{\sum_i p_i \teta_i^2}{1-\lmax(\PP)}.
  \label{eq:bnzf}
\end{align}
\begin{remark}
In Section \ref{sec:results}, our simulation results will illustrate that this upper bound may be several orders of magnitude too high than the true performance of the algorithm.
\end{remark}


\section{Average over the noise}
\label{sec:gauss}

Our exact MSE expression given in Proposition~\ref{prop:MSE} depends on the noise vector $\veta$. In practice, of course, $\veta$ is unknown, but we may have information about its statistics. In this section, we suppose that $\veta$ is drawn from a probability distribution: in particular, we assume that its elements $\eta_i$ are i.i.d. random variables with zero-mean and variance $\sigma^2$. Here, it is important to differentiate between two sources of randomness: the random row-selections made by the algorithm and the random vector $\veta$. In what follows, $\E$ is understood as the conditional expectation operator over the randomness of the algorithm, with $\veta$ fixed, and we define $\En$ as the average over the noise.

It is convenient to 
rewrite \eqref{eq:fixed_noise_iteration_MSE} as
\begin{equation}\label{eq:mse-s}
\begin{aligned}
  \MSE{k}&= 
  \opvec(\mI_{n})^T 
  \left[
  \QQ^k \left( \vzs{0}-\vv_{1} \right)
  \right.\\
  &\;\;
  \left.
  \quad+  
  f_k (\DD) \left( \vz^{(0)}-\vv_{2} \right)
  \right]
  +\tr \opmat(\vv_1)
  \;,
\end{aligned}
\end{equation}
where 
\begin{align}
  f_k(\DD)=\sum_{0 \le \ell < k} \QQ^\ell \DD \PP^{k-1-\ell}.
  \label{eq:fk}
\end{align}
Since $f_k(\DD)$ is a linear function, we have $\En  f_k(\DD) = f_k(\En \DD) = 0$.
Averaging \eqref{eq:mse-s} over the noise, we get the following proposition.
\begin{figure*}
\centering
\subfigure[]{
       \label{fig:gaussian_matrix1_results}
         \input{figs/fig_gaussian_matrix1_results.tikz}
        }   
        \hspace{2em}
\subfigure[]{
        \label{fig:tomography_matrix1_results}
          \input{figs/fig_tomography_matrix1_results.tikz}
     }
     \caption{(a)  The mean squared error $\E||\vx^{(k)} - \vx||^2$ is shown on a logarithmic scale as a function of the iteration number $k$.
            The matrix $\mA$ has Gaussian entries,
             and the error vector $\eta$ is fixed in advance, with $||\eta||^2 = 1.6$.
             The average results from 1007 trials are shown as the blue curve, and the results from 150 of the trials are
             shown in gray. The analytical expression 
             (\ref{eq:fixed_noise_iteration_MSE}) is shown as a dashed green line, and clearly matches the simulation results quite well.
             The Needell \cite{needell_paved_2014} and Zouzias-Freris \cite{zouzias_randomized_2013} bounds are shown as well, and are far higher
             than the true MSE.
             (b) The mean square error $\En  \E ||\vx^{(k)} - \vx||^2$ averaged over both the algorithm's randomness and the noise is 
             shown on a logarithmic scale as a function of the iteration number $k$. The matrix $\mA$ is the measurement matrix of a tomographic
             system (generated by the AIR Tools package \cite{hansen_air_2012}), and the error vector $\veta$ is a zero mean Gaussian vector with variance
               $2.25 \times 10^{-4}$, drawn independently with each trial. The average of 1007 trials
               are shown in blue along with the results from 150 of the trials in gray. The analytical expression for
               the Gaussian noise case (\ref{eq:mse-avg-noise}) clearly matches the simulation results. The noise-averaged Zouzias-Freris bound 
               is shown as well for comparison.
         }
    \label{fig:results}
\end{figure*}
\begin{proposition}
The MSE of RKA at the $k$th iteration averaged over both the randomness of the algorithm and noise is
\begin{equation}\label{eq:mse-avg-noise}
\begin{aligned}
  \En  \MSE{k} 
  &=
  \opvec(\mI_{n})^T 
  \left[
    \QQ^k \left( \vzs{0} - \En  \vv_{1} \right)
  \right.  \\
  &\;\;\;\;
  \left.
  \quad\quad- \En f_k (\DD) \vv_{2}
  \right]   +\tr \opmat(\En \vv_1).
\end{aligned}
\end{equation}
\end{proposition}
This formula involves two noise-related quantities, 
$\En \vv_1$ and
 $\En f_k (\DD) \vv_2 $,
 both of which are second-order in the noise. This shows that our knowledge of the second-order statistics of the noise is sufficient to compute them.
In particular, the first term is given by 
\[
  \En  \vv_1 = \sigma^2 (\mI - \QQ) ^{-1} 
  \left[
    \sum_{i}p_i 
    \frac{\tva_i^{ \otimes 2}}{||\va_i||^2}
    + g\left( (\mI-\mP)^{-1} \right)
      \right],
\]
where we define the matrix function
\begin{equation}
    g(\mM) = \sum_i \frac{p_i^2}{||\va_i||^2} \left( \tva_i \otimes \projp{i} + \projp{i} \otimes \tva_i \right) \mM \tva_i.
\end{equation}
(In these expressions the extra factors of $||\va_i||^2$ are not erroneous---they account for the varying signal-to-noise ratio of the measurements.) 

The second noise-related term is computed by
\[
\begin{aligned}
  \En  f_k (\DD) \vv_{2}
  &=\sum_{0 \le \ell < k} \QQ^\ell \hairspace \En  \DD \PP^{k-1-\ell} (\mI - \PP)^{-1} \vf
  \\
  &= \sigma^2 \sum_{0 \le \ell < k} \QQ^\ell \hairspace \En  g( \PP^{k-1-\ell} (\mI - \PP)^{-1}).
\end{aligned}
\]
\begin{remark}
  The first term on the right-hand side of \eqref{eq:mse-avg-noise} decays exponentially because $\lambda_{\max}(\mQ) < 1$. Thus, the limiting MSE averaged over both the randomness of the algorithm and noise is $\En  \MSE{\infty} = \tr \opmat(\En  \vv_1)$.
\end{remark}

\section{Experimental Results}
\label{sec:results}

We verified our results with numerical simulations. We took care in our implementations of the matrices $\QQ$ and $\DD$ in order to minimize
the time- and space-complexity. $\QQ$ is an $n^2 \times n^2$ matrix, which in a naive implementation would require $O(n^4)$ storage and 
$O(n^4)$ computation to multiply by a vector. Instead, we use the fact that
\begin{equation}
    \QQ \vx = \opvec \left(\sum_i p_i \projp{i} \opmat(\vx) \projp{i} \right)
\end{equation}
to implement multiplication by $\QQ$ with no additional storage in time $O(m n^2)$ [since $\projp{i}$ can be
multiplied by other matrices in time $O(n^2)]$. Meanwhile, we use the fact that
\begin{equation}
    \DD \vx = \opvec \left( \sum_i \frac{p_i \veta_i}{||\va_i||^2} (\projp{i}\vx\va_i^T  + \va_i \vx^T \projp{i})   \right)
\end{equation}
to implement $\DD$ without any additional storage. This saves no computation, since it takes $O(m n^2)$ time.

For the noise-averaged formula \eqref{eq:mse-avg-noise},
we can use the structure of $\QQ$ to compute the complex term $\En  f_k (\DD) \vv_{2}$ in $O(kmn^2)$ time.
Alternatively, we could use an eigenvector decomposition of $\mP$ and $\mQ$ to compute
it in a time constant in $k$: we must use $O(n^4)$ space and $O(m n^4)$ time. This would
make sense if we wanted to compute the MSE for a single, moderately large $k$.

The results of two experiments are shown in this paper. 
First, we tested the fixed noise formula (\ref{eq:fixed_noise_iteration_MSE}). We drew a single noise vector $\veta$ with $||\veta||^2 = 1.6$,
and a starting error $\vz^{(0)}$,
and choose a $150 \times 50$ measurement matrix $\mA$ that had i.i.d. Gaussian entries. Then we ran 1007 separate
trials of the randomized Kaczmarz algorithm, with each trial running for 2000 iterations and starting with an error vector $\vz^{(0)}$.
We plotted the average MSE of the trials at each iteration on a log scale. The results are shown in Figure \ref{fig:gaussian_matrix1_results},
and show that the expression we derived in (\ref{eq:fixed_noise_iteration_MSE}) matches the numerical results very well.
We also plotted existing bounds \cite{needell_randomized_2010,zouzias_randomized_2013}
as well. The bounds are significantly higher than the true MSE.

Next, we tested the noise-averaged formula (\ref{eq:mse-avg-noise}). We used the AIR Tools package in MATLAB \cite{hansen_air_2012} to
generate a tomography measurement matrix $\mA$ of size $148 \times 100$. The noise vector $\veta$ had i.i.d.~entries with variance
$\sigma^2 = 2.25\times10^{-4}$ and was drawn independently for each trial. We ran 1007 separate trials of the randomized Kaczmarz algorithm,
with each trial running for 3000 iterations. The results are shown in Figure \ref{fig:tomography_matrix1_results}. The close match between empirical
and theoretical curves verify our expression for the noise-averaged MSE (\ref{eq:mse-avg-noise}). The graph also shows that
the noise-averaged version of the Zouzias-Freris bound is more than two orders of magnitude higher than the true limiting MSE in this case.

\section{Conclusions}
\label{sec:conclusion}
We provided a complete characterization of the randomized Kaczmarz algorithm when applied to inconsistent linear systems. We developed an exact
formula for the MSE of the algorithm when the measurement vector is corrupted by a fixed noise vector. We also showed how to average this expression
over a noise distribution with known first and second-order moments. We described efficient numerical implementations of these expressions that limit the time- and space-complexity. Simulations show that the exact MSE expressions we derived have excellent matches with the numerical results. Moreover, our experiments indicate that existing upper bounds on the MSE may be loose by several orders of magnitude.


\providecommand{\url}[1]{}
\renewcommand{\url}[1]{}
\bibliographystyle{IEEEtran}


\end{document}